\magnification 1200
\input plainenc
\input amssym
\fontencoding{OT1}
\inputencoding{latin3}

\tolerance 7500
\relpenalty 10000
\binoppenalty 10000
\parindent 1.5em

\hsize 17.25truecm
\vsize 24.5truecm
\hoffset 0truecm
\voffset -1.5truecm

\font\TITLE cmbx10 at 14.4pt
\font\tenrm cmr10
\font\cmtenrm cmr10
\font\tenit cmti10
\font\tenbf cmbx10
\font\teni cmmi10 \skewchar\teni '177
\font\tensy cmsy10 \skewchar\tensy '60
\font\tenex cmex10
\font\teneufm eufm10
\font\eightrm cmr8
\font\cmeightrm cmr8
\font\eightit cmti8
\font\eightbf cmbx8
\font\eighti cmmi8 \skewchar\eighti '177
\font\eightsy cmsy8 \skewchar\eightsy '60
\font\eightex cmex8
\font\eighteufm eufm8
\font\cmsevenrm cmr7

\font\cmsixrm cmr6

\font\sixbf cmbx6
\font\sixi cmmi6 \skewchar\sixi '177
\font\sixsy cmsy6 \skewchar\sixsy '60
\font\sixeufm eufm6

\font\cmfiverm cmr5

\font\fivebf cmbx5
\font\fivei cmmi5 \skewchar\fivei '177
\font\fivesy cmsy5 \skewchar\fivesy '60
\font\fiveeufm eufm5
\font\tencmmib cmmib10 \skewchar\tencmmib '177
\font\eightcmmib cmmib8 \skewchar\eightcmmib '177
\font\sevencmmib cmmib7 \skewchar\sevencmmib '177
\font\sixcmmib cmmib6 \skewchar\sixcmmib '177
\font\fivecmmib cmmib5 \skewchar\fivecmmib '177
\font\kirrm wncyr8
\font\kirit wncyi8

\newfam\cmmibfam
\textfont\cmmibfam\tencmmib \scriptfont\cmmibfam\sevencmmib
\scriptscriptfont\cmmibfam\fivecmmib
\def\tenpoint{\def\rm{\fam0\tenrm}\def\it{\fam\itfam\tenit}%
	\def\bf{\fam\bffam\tenbf}
	\textfont0\cmtenrm \scriptfont0\cmsevenrm \scriptscriptfont0\cmfiverm
  	\textfont1\teni \scriptfont1\seveni \scriptscriptfont1\fivei
  	\textfont2\tensy \scriptfont2\sevensy \scriptscriptfont2\fivesy
  	\textfont3\tenex \scriptfont3\tenex \scriptscriptfont3\tenex
  	\textfont\itfam\tenit
	\textfont\bffam\tenbf \scriptfont\bffam\sevenbf
	\scriptscriptfont\bffam\fivebf
	\textfont\eufmfam\teneufm \scriptfont\eufmfam\seveneufm
	\scriptscriptfont\eufmfam\fiveeufm
	\textfont\cmmibfam\tencmmib \scriptfont\cmmibfam\sevencmmib
	\scriptscriptfont\cmmibfam\fivecmmib
	\normalbaselineskip 12pt
	\setbox\strutbox\hbox{\vrule height8.5pt depth3.5pt width0pt}%
	\normalbaselines\rm}
\def\eightpoint{\def\rm{\fam 0\eightrm}\def\it{\fam\itfam\eightit}%
	\def\bf{\fam\bffam\eightbf}%
	\textfont0\cmeightrm \scriptfont0\cmsixrm \scriptscriptfont0\cmfiverm
	\textfont1\eighti \scriptfont1\sixi \scriptscriptfont1\fivei
	\textfont2\eightsy \scriptfont2\sixsy \scriptscriptfont2\fivesy
	\textfont3\eightex \scriptfont3\eightex \scriptscriptfont3\eightex
	\textfont\itfam\eightit
	\textfont\bffam\eightbf \scriptfont\bffam\sixbf
	\scriptscriptfont\bffam\fivebf
	\textfont\eufmfam\eighteufm \scriptfont\eufmfam\sixeufm
	\scriptscriptfont\eufmfam\fiveeufm
	\textfont\cmmibfam\eightcmmib \scriptfont\cmmibfam\sixcmmib
	\scriptscriptfont\cmmibfam\fivecmmib
	\normalbaselineskip 11pt
	\abovedisplayskip 5pt
	\belowdisplayskip 5pt
	\setbox\strutbox\hbox{\vrule height7pt depth2pt width0pt}%
	\normalbaselines\rm
}

\def\empty{}

\catcode`\@ 11
\catcode`\" 12

\def\hat{\mathaccent"705E }

\def\newl@bel#1#2{\expandafter\def\csname l@#1\endcsname{#2}}
\openin 11\jobname .aux
\ifeof 11
	\closein 11\relax
\else
	\closein 11
	\input \jobname .aux
	\relax
\fi

\newcount\c@section
\newcount\c@subsection
\newcount\c@subsubsection
\newcount\c@equation
\newcount\c@bibl
\c@section 0
\c@subsection 0
\c@subsubsection 0
\c@equation 0
\c@bibl 0
\def\lab@l{}
\def\label#1{\immediate\write 11{\string\newl@bel{#1}{\lab@l}}%
	\ifhmode\unskip\fi}

\def\section#1{\global\advance\c@section 1
	{\par\vskip 3ex plus 0.5ex minus 0.1ex
	\rightskip 0pt plus 1fill\leftskip 0pt plus 1fill\noindent
	{\bf\S\thinspace\number\c@section .~#1}\par\penalty 25000%
	\vskip 1ex plus 0.25ex}
	\gdef\lab@l{\number\c@section.}
	\c@subsection 0
	\c@subsubsection 0
	\c@equation 0
}
\def\subsection{\global\advance\c@subsection 1
	\par\vskip 1ex plus 0.1ex minus 0.05ex{\bf\number\c@subsection. }%
	\gdef\lab@l{\number\c@section.\number\c@subsection}%
	\c@subsubsection 0\c@equation 0%
}
\def\subsubsection{\global\advance\c@subsubsection 1
	\par\vskip 1ex plus 0.1ex minus 0.05ex%
	{\bf\number\c@subsection.\number\c@subsubsection. }%
	\gdef\lab@l{\number\c@section.\number\c@subsection.%
		\number\c@subsubsection}%
}
\def\equation{\global\advance\c@equation 1
	\gdef\lab@l{\number\c@section.\number\c@subsection.%
	\number\c@equation}{\rm\number\c@equation}
}
\def\bibitem#1{\global\advance\c@bibl 1
	[\number\c@bibl]%
	\gdef\lab@l{\number\c@bibl}\label{#1}
}
\def\ref@ref#1.#2:{\def\REF@{#2}\ifx\REF@\empty{\S\thinspace#1}%
	\else\ifnum #1=\c@section {#2}\else {\S\thinspace#1.#2}\fi\fi
}
\def\ref@eqref#1.#2.#3:{\ifnum #1=\c@section\ifnum #2=\c@subsection
	{(#3)}\else{#2\thinspace(#3)}\fi\else{\S\thinspace#1.#2\thinspace(#3)}\fi
}
\def\ref#1{\expandafter\ifx\csname l@#1\endcsname\relax
	{\bf ??}\else\edef\mur@{\csname l@#1\endcsname :}%
	{\expandafter\ref@ref\mur@}\fi
}
\def\eqref#1{\expandafter\ifx\csname l@#1\endcsname\relax
	{(\bf ??)}\else\edef\mur@{\csname l@#1\endcsname :}%
	{\expandafter\ref@eqref\mur@}\fi
}
\def\cite#1{\expandafter\ifx\csname l@#1\endcsname\relax
	{\bf ??}\else\hbox{\bf\csname l@#1\endcsname}\fi
}

\def\land{\mathbin{\&}}
\def\dom{\mathop{\rm dom}}
\def\mes{\mathop{\rm mes}}
\catcode`\@ 12

\def\proof{\par\medskip{\rm P$\,$r$\,$u$\,$v$\,$o.}\ }
\def\endproof{{\parfillskip 0pt\hfill$\square$\par}\medskip}

\immediate\openout 11\jobname.aux


\frenchspacing\tenpoint
\leftline{\ }\vskip 0.25cm
{\leftskip 0cm plus 1fill\rightskip 0cm plus 1fill\parindent 0cm\baselineskip 15pt
\TITLE Pri komparo de integraloj de Lebesgue kaj Riemann\par
en konstrua matematika analizo
\par\vskip 0.15cm\rm A.$\,$A.~Vladimirov\par}
\vskip 0.25cm
$$
	\vbox{\hsize 0.75\hsize\leftskip 0cm\rightskip 0cm
	\eightpoint\rm
	{\bf Resumo:\/} Estas trovita sumebleco de ĉiu konstrua funkcio, kiu estas
	difinita preskaŭ ĉie en la intervalo $[0,1]$ kaj integralebla en senco
	de Riemann.\par
	}
$$

\vskip 0.5cm
\section{Enkonduko}\label{par:1}%
\subsection
Teorio de integralo, disvolvata enkadre de konstrua matematika analizo, en multaj
flankoj okazas esence diferenca de tradicia teorio, kiu apogas sin sur la teorio de aroj.
Ekzemple, tie disiras [\cite{Vl:2012}] du sencoj de integraleblo, ligitaj al la nomoj
de Riemann kaj Darboux, kiuj laŭ la vidpunkto de teorio de aroj estas rigardataj
kiel tute interŝanĝeblaj. Sekve, la problemo pri komparo de integraloj de Riemann
kaj Lebesgue, en tradicia teorio facile solvata kun helpo de la fakto de sumebleco
de ĉiu funkcio, integralebla en senco de Darboux, laŭ la vidpunkto de konstrua
branĉo de la matematiko okazas sufiĉe neklara.

Uzante rezultaton de G.$\,$S.~Cejtin pri senrompeco de konstruaj funkcioj
[\cite{Ku:1973}: Ĉap.~9, \S$\,$2], oni povas nemalfacile trovi, ke kiu ajn
{\it ĉie\/} difinita konstrua funkcio $f\colon[0,1]\to\Bbb R$ en okazo de ĝia
integraleblo en senco de Riemann estas sumebla. Tamen tiu rezultato ne elĉerpas
la problemon. Enkadre de teorio de integralo, multe pli nature estas konsideri
ne funkciojn, difinitajn ĉie (kiel estis farate, ekzemple, en la verko
[\cite{Dem:1967}]), sed funkciojn, difinitajn {\it preskaŭ ĉie\/} (precizaj difinoj
estos donataj poste). En okazo de konsiderado de tiaj funkcioj, antaŭcitita pruvo
ne povas esti reproduktata senpere. Tiu ĉi artikolo estas dediĉita al la trovo
de sumebleco de ĉiu funkcio, kiu estas difinita preskaŭ ĉie en la intervalo $[0,1]$
kaj integralebla en senco de Riemann.

\subsection
Poste ni apogos nin sur konstrua formo de tiu difinmaniero de integralo de Lebesgue,
kiu estis ellaborita de P.~Daniell [\cite{RN:1972}: \S$\,$61]. Tiu formo estas
preskaŭ egala kun "klasika", postulante nur precizigon de kelke da esprimoj
pro plenumata en la konstrua analizo akcepto de ekzistado de Specker'aj nombraj
vicoj. La skiza esprimo de tiu teorio estos donata en~\ref{par:2}.


\section{Funkcioj, difinitaj preskaŭ ĉie}\label{par:2}
\subsection
Vicon $\{h_n\}_{n=0}^\infty$ da nemalpozitivaj multangulaj funkcioj [NMF]
ni poste nomos {\it regula\/} [RVNMF], se ĝi obeas rilaton
$$
	(\forall n\in\Bbb N)\qquad\int_0^1 h_n<2^{-n}.
$$
Fundamento de la teorio de integralo de Lebesgue konsistas el la sekvanta
facila fakto.

\subsubsection\label{prop:1.1}
{\it Estu $h$ multangula funkcio, kaj estu $\{h_n\}_{n=0}^\infty$ {\rm RVNMF},
kiu obeas la kondiĉon
$$
	\int_0^1 h>\sum_{n=0}^\infty\int_0^1 h_n.
$$
En tiu okazo estas realigebla punkto $\xi\in [0,1]$, por kiu nombra serio
$\sum_{n=0}^\infty h_n(\xi)$ celas al iu reala nombro, severe malgranda
ol la nombro $h(\xi)$.
}%
\proof
Estas klare, ke estas realigebla valoro $\varepsilon\in (0,1)$, kiu obeas
la kondiĉon
$$
	\int_0^1 (h-\varepsilon)>\sum_{n=0}^\infty\int_0^1(1+\varepsilon)^nh_n.
$$
Do, per ordinara maniero de sekcado de intervalo, oni povas facile trovi
realigeblecon de punkto $\xi\in [0,1]$ kun eco
$$
	(\forall m\in\Bbb N)\qquad h(\xi)-\varepsilon\geqslant
		\sum_{n=0}^m(1+\varepsilon)^nh_n(\xi).
$$
Tiu punkto ja estas serĉata.
\endproof

\subsection
Kun ĉiu RVNMF $\frak h\rightleftharpoons\{h_n\}_{n=0}^\infty$ oni povas ligi
aron $M_{\frak h}$ da duopoj $(x,\gamma)\in [0,1]\times\Bbb R$ obeantaj al
la kondiĉo
$$
	(\forall m\in\Bbb N)\qquad\sum_{n=0}^m h_n(x)\leqslant\gamma.
$$
Simbolo $\Omega_{\frak h}$ signifos poste ombron $\{x\in [0,1]\;:
\;(\exists\gamma\in\Bbb R)\;(x,\gamma)\in M_{\frak h}\}$ de la aro $M_{\frak h}$.
La ideo pri {\it preskaŭ plena\/} subaro [PPS] de la intervalo $[0,1]$ povas nun
esti esprimata per la sekvanta maniero.

\subsubsection\label{prop:2.1}
{\it Subaro $X\subseteq [0,1]$ estas nomata preskaŭ plena, se oni povas realigi
je \hbox{\rm RVNMF} $\frak h$ kun eco $\Omega_{\frak h}\subseteq X$.
}

\medskip
Simile al la situacio de "klasika" analizo, oni povas facile konvinkiĝi pri
la sekvantaj faktoj.

\subsubsection
{\it Ĉiu {\rm PPS} enhavas almenaŭ unu punkton} [\ref{prop:1.1}].

\subsubsection\label{prop:2.2.3}
{\it Komuna parto de du arbitraj {\rm PPS} ankaŭ estas {\rm PPS}.
}

\subsubsection\label{prop:2.2.5}
{\it Komuna parto
$
	\bigcap_{n=0}^\infty X_n=\{x\in [0,1]\;:\;(\forall n\in\Bbb N)\quad (n,x)\in X\}
$
de arbitra vico $X\subseteq\Bbb N\times [0,1]$ da {\rm PPS}
$
	X_n\rightleftharpoons\{x\in [0,1]\;:\; (n,x)\in X\}
$
ankaŭ estas {\rm PPS}.
}%
\proof
Estu $\{h_{n,k}\}_{k=0}^\infty$ RVNMF, kiuj en senco de la difino~\ref{prop:2.1}
respondas al aroj $X_n$. Ni konsideru vicon $\frak g\rightleftharpoons
\{g_k\}_{k=0}^\infty$ da NMF kun formo
$$
	g_k\rightleftharpoons\sum_{n=0}^k 2^{-2n-1}\,h_{n,k-n}.
$$
Tiu ĉi vico videble estas regula. Dume ĉiu duopo $(x,\gamma)\in M_{\frak g}$
obeas rilatojn
$$
	\leqalignno{\sum_{k=0}^m h_{n,k}(x)&\leqslant
		2^{2n+1}\,\sum_{k=0}^m g_{n+k}(x)\cr &\leqslant 2^{2n+1}\gamma,}
$$
kiuj garantias verecon de rilatoj $x\in X_n$. Sekve, la rilato
$\Omega_{\frak g}\subseteq\bigcap_{n=0}^\infty X_n$ estas ankaŭ vera.
\endproof

\subsubsection\label{prop:jihe}
{\it Por ĉiu {\rm RVNMF} $\{h_n\}_{n=0}^\infty$ oni povas realigi je {\rm PPS\/}
$X\subseteq [0,1]$ kun eco
$$
	(\forall x\in X)\;(\exists C\in\Bbb R)\;(\forall n\in\Bbb N)\qquad
		h_n(x)<C\cdot(3/4)^n.
$$
}%
\proof
Ni konsideru vicon $\frak g\rightleftharpoons\{g_n\}_{n=0}^\infty$ da {\rm NMF}
kun formo
$$
	g_n\rightleftharpoons{(4/3)^{2n}h_{2n}+(4/3)^{2n+1}h_{2n+1}\over 2}.
		\leqno(\equation)
$$\label{eq:gn}%
Tiu ĉi vico estas regula pro la rilatoj
$$
	\eqalign{\int_0^1 g_n&<{(4/3)^{2n}\,2^{-2n}+(4/3)^{2n+1}\,2^{-2n-1}
			\over 2}\cr &=(5/6)\cdot (4/9)^n\cr &\leqslant 2^{-n},}
$$
garantiataj de reguleco de la vico $\{h_n\}_{n=0}^\infty$. Dume ĉiu duopo
$(x,\gamma)\in M_{\frak g}$ obeas [\eqref{eq:gn}] kondiĉojn $h_n(x)\leqslant
2\gamma\cdot (3/4)^n$.
\endproof

\subsection
Ni konsideru tuj la situacion, se al ĉiu punkto $x\in X$ de iu PPS, ligita kun
la vico $\frak h$ en senco de la difino~\ref{prop:2.1}, respondas iu konstrua objekto
$f(x)$. Laŭ la vidpunkto de konstrua kompreno de matematikaj juĝoj, tio signifas,
ke estas realigebla kalkulplano $\hat f$, kiu konvertas ĉiun duopon $(x,\gamma)\in
M_{\frak h}$ al objekto $f(x)$. Tiu rimarko kondukas nin al la sekvanta esprimo
de la ideo pri {\it difinita preskaŭ ĉie\/} konstrua funkcio de reala argumento
[DPĈKF].

\subsubsection
{\it Difinita preskaŭ ĉie en la intervalo $[0,1]$ realvalora funkcio estas duopo
$\{\frak h,\hat f\}$, konsistanta el {\rm RVNMF} $\frak h$ kaj konstrua bildigo
$\hat f\colon M_{\frak h}\to\Bbb R$ kun eco
$$
	(\forall (x_1,\gamma_1),\,(x_2,\gamma_2)\in M_{\frak h})\qquad
		\bigl(\hat f(x_1,\gamma_1)\neq\hat f(x_2,\gamma_2)\bigr)\to
		\bigl(x_1\neq x_2\bigr).
$$
}

Per simbolo $f(x)$, kie $f=\{\frak h,\hat f\}$ kaj $x\in\Omega_{\frak h}$, ni poste
signos valoron $\hat f(x,\gamma)$, respondantan al arbitre fiksita valoro
$\gamma\in\Bbb R$ kun eco $(x,\gamma)\in M_{\frak h}$. En okazo, se RVNMF $\frak h$
estas ligata kun iu DPĈKF $f=\{\frak h,\hat f\}$, la aro $\Omega_{\frak h}$ estos
ankaŭ signata per simbolo $\dom f$.

\subsection
La ideo pri {\it sumebleco\/} de DPĈKF povas esti enkondukita per la sekvanta
maniero.

\subsubsection\label{prop:sum}
{\it {\rm DPĈKF} $f$ estas nomata sumebla, se estas realigebla vico da multangulaj
funkcioj $\{f_n\}_{n=0}^\infty$ kun eco
$$
	(\forall n\in\Bbb N)\qquad \int_0^1 |f_{n+1}-f_n|<2^{-n},\leqno(\equation)
$$\label{eq:sum0}%
por kiu sur iu {\rm PPS} estas vera la egaleco
$$
	f(x)\equiv\lim_{n\to\infty} f_n(x).\leqno(\equation)
$$\label{eq:sum1}%
Integralo de Lebesgue de {\rm DPĈKF} $f$ estas en tiu okazo la nomo de la valoro
$$
	\int_0^1 f\rightleftharpoons\lim_{n\to\infty}\int_0^1 f_n.
$$
}%

\penalty -5000
Oni povas trovi sekvantajn kvin faktojn.

\subsubsection\label{prop:3.2}
{\it Estu {\rm DPĈKF} $f$ sumebla kaj havanta severe pozitivan integralon.
Do estas realigebla punkto $\xi\in\dom f$ kun eco $f(\xi)>0$.
}%
\proof
Ni fiksu je RVNMF $\frak h=\{h_n\}_{n=0}^\infty$, por kiu sur la aro
$\Omega_{\frak h}\subseteq\dom f$ estas vera la egaleco~\eqref{eq:sum1}.
Konforme al faritaj supozoj, estas realigebla nombro $m\geqslant 1$ kun eco
$$
	\int_0^1 f_m>\sum_{k=0}^\infty\int_0^1
		\Bigl[|f_{m+k+1}-f_{m+k}|+2^{-m}h_k\Bigr].
$$
Realigebleco de punkto $\xi\in\Omega_{\frak h}$ kun eco serĉata estas nun
garantiata de la juĝo~\ref{prop:1.1}.
\endproof

\subsubsection
{\it Valoro de integralo de Lebesgue de ĉiu sumebla {\rm DPĈKF} estas difinita
unusence} [\ref{prop:sum}, \ref{prop:2.2.3}, \ref{prop:3.2}].

\subsubsection\label{prop:lev}
{\it Estu vico $\{f_n\}_{n=0}^\infty$ da sumeblaj {\rm DPĈKF} obeanta al la
rilato~{\rm\eqref{eq:sum0}\/}. Do estas realigebla sumebla {\rm DPĈKF} $f$,
obeanta al la rilatoj
$$
	\displaylines{(\forall x\in\dom f)\qquad f(x)=\lim_{n\to\infty} f_n(x),\cr
		\int_0^1 f=\lim_{n\to\infty}\int_0^1 f_n.}
$$
}%
\proof
Ni ligu kun sumeblaj funkcioj $f_n$ esprimojn
$$
	(\forall x\in X)\qquad f_n(x)=\lim_{k\to\infty} f_{n,k}(x),
	\leqno(\equation)
$$\label{eq:lev0}%
kiuj estu veraj sur iu fiksita [\ref{prop:2.2.5}] PPS $X\subseteq [0,1]$,
kaj kies multangulaj funkcioj $f_{n,k}$ obeu kondiĉojn
$$
	\int_0^1 |f_{n,k+1}-f_{n,k}|<2^{-k}.\leqno(\equation)
$$\label{eq:lev1}%
Ni ankaŭ enkonduku en la konsideradon vicon $\frak g=\{g_n\}_{n=0}^\infty$
da NMF kun formo
$$
	g_n\rightleftharpoons\sum_{k=0}^{n+2} |f_{n+2-k,n+4+k}-f_{n+2-k,n+2+k}|.
	\leqno(\equation)
$$\label{eq:lev2}%
Laŭ garantiataj de la kondiĉoj~\eqref{eq:lev1} rilatoj
$$
	\leqalignno{\int_0^1 g_n&\leqslant\sum_{k=0}^{n+2}\,
		[2^{-n-3-k}+2^{-n-2-k}]\cr &<2^{-n},
	}
$$
la vico $\frak g$ estas regula. Tio signifas [\ref{prop:jihe}, \ref{prop:2.2.3},
\eqref{eq:lev2}], ke estas realigebla PPS $Y\subseteq X$ kun eco
$$
	(\forall x\in Y)\;(\exists C\in\Bbb R)\;(\forall n,k\in\Bbb N)\qquad
		|f_{n+2,n+4+2k}(x)-f_{n+2,n+2+2k}(x)|<C\cdot(3/4)^{n+k},
$$
kaj tial [\eqref{eq:lev0}] kun ankaŭ eco
$$
	(\forall x\in Y)\;(\exists C\in\Bbb R)\;(\forall n\in\Bbb N)\qquad
		|f_{n+2}(x)-f_{n+2,n+2}(x)|<4C\cdot(3/4)^n.\leqno(\equation)
$$\label{eq:lev5}%
Aliflanke, la kondiĉoj
$$
	\leqalignno{\int_0^1|f_{n+3,n+3}-f_{n+2,n+2}|&\leqslant
		\int_0^1\Bigl[|f_{n+3,n+3}-f_{n+3}|+|f_{n+3}-f_{n+2}|+
			|f_{n+2}-f_{n+2,n+2}|\Bigr]\cr
		&<2^{-n-2}+2^{-n-2}+2^{-n-1}\cr
		&=2^{-n}\cr}
$$
garantias [\ref{prop:jihe}] celhavecon de RVNMF $\{f_{n+2,n+2}\}_{n=0}^\infty$
sur iu PPS $Z\subseteq [0,1]$. Arbitra DPĈKF $f$ kun ecoj $\dom f\subseteq
Y\cap Z$ kaj
$$
	(\forall x\in\dom f)\qquad f(x)=\lim_{n\to\infty}f_{n+2,n+2}(x)
$$
ja estas serĉata [\eqref{eq:lev5}].
\endproof

\subsubsection\label{prop:nulint}
{\it Ĉiu nemalpozitiva sumebla {\rm DPĈKF} $f$ kun egala al $0$ integralo akiras
valoron $0$ preskaŭ ĉie.
}%
\proof
Estas sufiĉe aluzi la juĝon~\ref{prop:lev} al la funkcia vico $\{nf\}_{n=0}^\infty$.
\endproof

\subsection
La ideon pri preskaŭ ĉie difinita {\it karakteriga funkcio\/} de subaro
ni enkondukas per la sekvanta maniero.

\subsubsection
{\it {\rm DPĈKF} $f\rightleftharpoons\{\frak h,\hat f\}$ estas nomata karakteriga
funkcio de subaro $X\subseteq [0,1]$, se ĉiu duopo $(x,\gamma)\in M_{\frak h}$
obeas rilaton
$$
	\biggl(\bigl(x\in X\bigr)\land\bigl(\hat f(x,\gamma)=1\bigr)\biggr)\lor
		\biggl(\bigl(x\not\in X\bigr)\land\bigl(\hat f(x,\gamma)=0\bigr)
		\biggr).
$$
}%

\medskip
Subaron $X\subseteq [0,1]$ ni nomas {\it mezurebla\/}, se por ĝi estas realigebla
sumebla karakteriga funkcio. Integralo de tiu karakteriga funkcio estas nomata
{\it mezuro de Lebesgue\/} de subaro $X$, kaj estas signata per simbolo $\mes X$.

\subsubsection\label{prop:aup}
{\it  Ĉiu mezurebla subaro kun nenula mezuro de Lebesgue enhavas almenaŭ
unu punkton} [\ref{prop:3.2}].

\subsubsection\label{prop:ppp}
{\it Subaro $X\subseteq [0,1]$ estas {\rm PPS} en tiu kaj nur en tiu okazo,
se ĝi estas mezurebla kaj ĝia mezuro estas egala al $1$} [\ref{prop:nulint}].


\section{Ĉefaj rezultatoj}\label{par:3}
\subsection
Ni fiksu arbitran DPĈKF $f=\{\{h_n\}_{n=0}^\infty,\hat f\}$, kaj ankaŭ iun ligatan
kun ĝi vicon $\{\Delta_n\}_{n=0}^\infty$ da mezureblaj aroj, obeantaj al rilatoj
$$
	(\forall x\in\Delta_n)\qquad h_n(x)<(2/3)^n
$$
kaj $\mes\Delta_n>1-(3/4)^n$. Per simboloj $\Gamma_n$ ni poste signos mezureblajn
[\ref{prop:lev}] arojn $\bigcap_{k=n}^\infty\Delta_k$, kiuj klare obeas kondiĉojn
$\mes\Gamma_n>1-4\cdot(3/4)^n$.

Ni fiksu tian solveblan aron $\Theta\subseteq\Bbb N^3$ da triopoj $(k,m,n)$
kun eco $k<2^m$, ke en la okazo $(k,m,n)\in\Theta$ intervalo $I_{k,m}
\rightleftharpoons (k\cdot 2^{-m},[k+1]\cdot 2^{-m})$ obeu kondiĉon
$\mes(I_{k,m}\cap\Gamma_n)>0$, kaj alie ĝi obeu kondiĉon $\mes(I_{k,m}\cap
\Gamma_n)<4^{-m}$. Ni fiksu [\ref{prop:aup}] ankaŭ kalkulplanon, kiu konvertas
ĉiun triopon $(k,m,n)\in\Theta$ al iu punkto $\zeta_{k,m,n}\in I_{k,m}\cap
\Gamma_n$, kaj konvertas ĉiun alian triopon $(k,m,n)\in\Bbb N^3$ kun eco $k<2^m$
al iu punkto $\zeta_{k,m,n}\in I_{k,m}\cap\dom f$.

Fine, ni fiksu direktitan aron $\frak A$, kies anoj havas formon
$$
	\alpha\rightleftharpoons(m,\,\{(k_l,m_l,n_l)\}_{l=0}^{2^m-1})
	\leqno(\equation)
$$\label{eq:frakA}%
kaj obeas kondiĉojn $m_l\geqslant m$ kaj $l\cdot 2^{-m}\leqslant
k_l\cdot 2^{-m_l}<(l+1)\cdot 2^{-m}$. Komparado de duopoj kun formo~\eqref{eq:frakA}
estas plenumata per komparado de iliaj unuaj elementoj. Kun tiu ĉi direktita
aro oni povas ligi ĝeneraligitan vicon $\{f_\alpha\}_{\alpha\in\frak A}$ da sumeblaj
DPĈKF kun formo
$$
	f_\alpha\rightleftharpoons\sum_{l=0}^{2^m-1}f(\zeta_{k_l,m_l,n_l})
		\cdot\chi_{l,m}.
$$
Simbolo $\chi_{l,m}$ signifas tie ĉi sumeblan karakterigan funkcion de antaŭe
priskribita intervalo $I_{l,m}$.

\subsubsection\label{prop:intRiem}
{\it En okazo de integraleblo de {\rm DPĈKF} $f$ en senco de Riemann,
ĝeneraligita vico $\{f_\alpha\}_{\alpha\in\frak A}$ havas celon en senco
de meza valoro.
}%
\proof
Ni fiksu arbitran valoron $\varepsilon>0$ kaj objekton $\alpha\in\frak A$ kun eco
$$
	(\forall\alpha',\alpha''>\alpha)\qquad
		\int_0^1 (f_{\alpha'}-f_{\alpha''})<\varepsilon/2.\leqno(\equation)
$$\label{eq:gxenvic}%
Fiksinte ankaŭ du arbitrajn duopojn
$$
	\alpha'\rightleftharpoons(m',\,\{(k'_\lambda,m'_\lambda,n'_\lambda)%
		\}_{\lambda=0}^{2^{m'}-1})>\alpha,\qquad
	\alpha''\rightleftharpoons(m'',\,\{(k''_\lambda,m''_\lambda,n''_\lambda)%
		\}_{\lambda=0}^{2^{m''}-1})>\alpha,
$$
ni povas realigi du ligatajn duopojn
$$
	\beta^\pm\rightleftharpoons(m,\,\{(k^\pm_l,m^\pm_l,n^\pm_l)%
		\}_{l=0}^{2^m-1}),
$$
kie $m=\min\{m',m''\}$ kaj
$$
	\leqalignno{f(\zeta_{k^-_l,m^-_l,n^-_l})-\varepsilon/4&<\inf\left\{
		\inf_{\lambda\;:\;l\leqslant\lambda\cdot 2^{m-m'}<l+1}
		f(\zeta_{k'_\lambda,m'_\lambda,n'_\lambda}),
		\inf_{\lambda\;:\;l\leqslant\lambda\cdot 2^{m-m''}<l+1}
		f(\zeta_{k''_\lambda,m''_\lambda,n''_\lambda})\right\},\cr
	f(\zeta_{k^+_l,m^+_l,n^+_l})+\varepsilon/4&>\sup\left\{
		\sup_{\lambda\;:\;l\leqslant\lambda\cdot 2^{m-m'}<l+1}
		f(\zeta_{k'_\lambda,m'_\lambda,n'_\lambda}),
		\sup_{\lambda\;:\;l\leqslant\lambda\cdot 2^{m-m''}<l+1}
		f(\zeta_{k''_\lambda,m''_\lambda,n''_\lambda})\right\}.}
$$
Estas klare, ke $\beta^\pm>\alpha$, kaj ke preskaŭ ĉie estas plenumataj la kondiĉoj
$$
	|f_{\alpha'}(x)-f_{\alpha''}(x)|<f_{\beta^+}(x)-f_{\beta^-}(x)
		+\varepsilon/2.
$$
Sekve [\eqref{eq:gxenvic}], estas vera ankaŭ la rilato
$$
	(\forall\alpha',\alpha''>\alpha)\qquad
		\int_0^1|f_{\alpha'}-f_{\alpha''}|<\varepsilon,
$$
kiu, laŭ la juĝo~\ref{prop:lev} kaj fakto de arbitreco de valoro $\varepsilon>0$,
garantias verecon de pruvata juĝo.
\endproof

\medskip
Ni rimarku, ke la formo de ideo pri konstrua integraleblo en senco de Riemann,
kiu estis uzata por la pruvo de la juĝo~\ref{prop:intRiem}, ne estas nure
ebla~[\cite{Ku:1970}: Ĉap.~3]. Tamen kvazaŭ pli liberigaj komprenoj efektive
kondukas al samsencaj ideoj~[\cite{Ku:1970}, \cite{Vl:2010}].

\subsection
En la daŭro de tiu ĉi subparagrafo ni supozas, ke konsiderata DPĈKF $f$ estas
integralebla en senco de Riemann. Per simbolo $g$ ni signas iun arbitre fiksitan
sumeblan celon de la ligita vico $\{f_\alpha\}_{\alpha\in\frak A}$ en senco
de meza valoro [\ref{prop:intRiem}]. Oni povas trovi la sekvantan fakton.

\subsubsection
{\it Por ĉiu $n\in\Bbb N$ oni povas realigi mezureblan aron $\Phi_n\subseteq
\dom f\cap\dom g$ kun eco $\mes\Phi_n>1-6\cdot(3/4)^n$, obeantan al la egaleco
$$
	(\forall x\in\Phi_n)\qquad f(x)=g(x).
$$
}%
\proof
Ni fiksu severe kreskantan vicon $\{m_k\}_{k=0}^\infty$ da nombroj $m_k\geqslant n$,
por kiu ligata vico $\{g_k\}_{k=0}^\infty$ da DPĈKF kun eco
$$
	g_k\rightleftharpoons\sum_{l=0}^{2^{m_k}-1}f(\zeta_{l,m_k,n})\cdot
		\chi_{l,m_k}
$$
obeas kondiĉon
$$
	(\forall k\in\Bbb N)\qquad\int_0^1|g_{k+1}-g_k|<2^{-k}.
$$
Krom tio, ni fiksu [\ref{prop:lev}, \ref{prop:nulint}] ankaŭ PPS $X\subseteq\dom g$,
kies ĉiu punkto obeas egalecon $\lim_{k\to\infty}g_k(x)=g(x)$.

Laŭ la konstruo de la aro $\Theta$, por ĉiuj $m\in\Bbb N$ estas plenumataj
la kondiĉoj
$$
	\mes\left(\Gamma_n\setminus\bigcup_{l\;:\;(l,m,n)\in\Theta}
		(I_{l,m}\cap\Gamma_n)\right)<2^{-m},
$$
kiuj signifas [\ref{prop:lev}], ke la aro
$$
	\Phi_n\rightleftharpoons X\cap\bigcap_{m=n}^\infty\left[
		\bigcup_{l\;:\;(l,m,n)\in\Theta}(I_{l,m}\cap\Gamma_n)\right]
$$
estas mezurebla kaj obeas rilaton $\mes\Phi_n>1-6\cdot(3/4)^n$. Dume, laŭ la konstruo
de la vico $\{g_k\}_{k=0}^\infty$, por ĉiu arbitre fiksita punkto $\xi\in\Phi_n$
estas realigebla celanta al ĝi vico $\{\xi_k\}_{k=0}^\infty$ da punktoj de la aro
$\Gamma_n$ kun eco
$$
	(\forall k\in\Bbb N)\qquad g_k(\xi)=f(\xi_k).
$$
Sekve, estas vera egaleco $\lim_{k\to\infty}f(\xi_k)=g(\xi)$. Aliflanke,
laŭ la konstruo de la aro $\Gamma_n$, estas realigebla valoro $\gamma\in\Bbb R$
kun eco
$$
	(\forall x\in\Gamma_n)\,(\forall m\in\Bbb N)\qquad
		\sum_{k=0}^m h_k(x)\leqslant\gamma.
$$
Tial, laŭ la juĝo pri senrompeco de konstruaj funkcioj [\cite{Ku:1973}:
Ĉap.~9, \S$\,$2], estas veraj ankaŭ egalecoj
$$
	\leqalignno{\lim_{k\to\infty} f(\xi_k)&=
		\lim_{k\to\infty}\hat f(\xi_k,\gamma)\cr
		&=\hat f(\xi,\gamma)\cr &=f(\xi).
	}
$$
Kolektante trovitajn rezultatojn, ni plenigas la pruvon.
\endproof

Por atingi la finan rezultaton pri sumebleco de konsiderata DPĈKF $f$,
ni tuj nur rimarku, ke la aro $\bigcup_{n=0}^\infty\Phi_n$, sur kiu DPĈKF $f$ estas
egala al sumebla DPĈKF $g$, estas preskaŭ plena [\ref{prop:ppp}].


\vskip 0.4cm
\eightpoint\rm
{\leftskip 0cm\rightskip 0cm plus 1fill\parindent 0cm
\bf Citita literaturo\par\penalty 20000}\vskip 0.4cm\penalty 20000
\bibitem{Vl:2012} {\kirit A.$\,$A.~Vladimirov\/}. {\kirrm O sravnenii integralov
Darbu i Rimana v konstruktivnom matematiqeskom analize~// Zapiski
nauq.~seminarov~POMI.~--- 2012.~--- T.~407.~--- S.~7--16.}

\bibitem{Ku:1973} {\kirit B.$\,$A.~Kuxner\/}. {\kirrm Lekcii po konstruktivnomu
matematiqeskomu analizu. M.: Nauka, 1973.}

\bibitem{Dem:1967} {\kirit O.~Demut.\/} {\kirrm Integral Lebega v konstruktivnom
analize~// Zapiski nauq.~seminarov~LOMI.~--- 1967.~--- T.~4.~--- S.~30--43.}

\bibitem{RN:1972} {\it F.~Riesz, B.~Sz.-Nagy\/}. Le\c cons d'analyse fonctionelle,
ed.~6. Budapest: Akad\' emiai Kiad\' o, 1972.

\bibitem{Ku:1970} {\kirit B.$\,$A.~Kuxner\/}. {\kirrm Nekotorye massovye problemy,
svyazannye s integrirovaniem konstruktivnyh funkcii~// Trudy MIAN
im.~V.$\,$A.~Steklova.~--- 1970.~--- T.~113.~--- S.~39--72.}

\bibitem{Vl:2010} {\kirit A.$\,$A.~Vladimirov\/}. {\kirrm O separabel{\char 126}nyh
napravlennostyah v konstruktivnyh topologiqeskih prostranstvah}~//
arXiv:1009.5335.~--- 2010.~--- {\kirrm S.~1--12.}
\bye